\documentclass[11pt]{article}
\usepackage{amsmath,amsthm,amssymb}
\usepackage[applemac]{inputenc}
\usepackage{mathrsfs}
\usepackage{graphicx,pict2e}
\usepackage{amsfonts}
\theoremstyle{plain}
\newtheorem{teo}{Theorem}[section]
\newtheorem{prop}[teo]{Proposition}

\theoremstyle{definition}
\newtheorem{ex}{Example}
\newtheorem{rem}{Remark}
\newtheorem{defi}{Definition}

\newtheorem{prob}{Problem}
\numberwithin{equation}{section}

\newcommand{\pical}{\mathcal{P}}
\newcommand{\tto}{\to}         
\newcommand{\res}{\mathop{\hbox{\vrule height 7pt width .5pt depth 0pt \vrule height .5pt width 6pt depth 0pt}}\nolimits}
\newcommand{\deb}{\rightharpoonup}

\newcommand{\haus}{\mathcal H}
\newcommand{\huno}{\haus^1}

\newcommand\io{\int_{\Omega}}

\newcommand{\M}{\mathcal M}

\newcommand{\lcal}{\mathcal L}

\newcommand{\W}{\mathcal W}

\newcommand{\R}{\mathbb R}

\DeclareMathOperator{\spt}{spt}
\DeclareMathOperator{\Lip}{Lip}

\DeclareMathOperator{\diam}{diam}

\title{Introduction to Optimal Transport Theory}
\author{Filippo Santambrogio\thanks{\scriptsize\ CEREMADE, UMR CNRS 7534, Universit\'e Paris-Dauphine, Pl. de Lattre de Tassigny, 75775 Paris Cedex 16, FRANCE
\texttt{filippo@ceremade.dauphine.fr , http://www.ceremade.dauphine.fr/$\sim$filippo}.}}
\date{Grenoble, June 15th, 2009\\Revised version : March 23rd, 2010}
\begin{document}
\maketitle

\begin{abstract}
These notes constitute a sort of Crash Course in Optimal Transport Theory. The different features of the problem of Monge-Kantorovitch are treated, starting from convex duality issues. The main properties of space of probability measures endowed with the distances $W_p$ induced by optimal transport are detailed. The key tools to put in relation optimal transport and PDEs are provided.
\end{abstract}



\tableofcontents

\section*{Introduction}

 These very short lecture notes do not want to be an exhaustive
presentation of the topic, but only a short list of results, concepts and ideas which are useful when dealing for the first time with the theory of Optimal Transport. Several of these ideas have been used, and explained in deeper details, during the other classes of the Summer School ``Optimal transportation : Theory and applications'' which were the occasion for the redaction of these notes. The style that was chosen when preparing them, in view of their use during the Summer School, was highly informal and this revised version will respect the same style. 

The main
references for the whole topic are the two books on the subject by C. Villani (\cite{villani,villani06}). For what concerns curves in the space of probability measures, the best specifically focused reference is \cite{AmbGigSav}. Moreover, I'm also very indebted to the approach that L. Ambrosio used in a course at SNS Pisa in 2001/02
and I want to cite as another possible reference \cite{dispensediambrosio}.

The motivation for the whole subject is the following problem proposed by Monge in 1781 (\cite{Monge}): given two densities of mass $f,\,g\geq 0$ on $\R^d$, with $\int f=\int g =1$, find a map $T:\R^d\to\R^d$ pushing the first one onto the other, i.e.
such that
\begin{equation}\label{push-forward}
\int_{A}g(x)dx=\int_{T^{-1}(A)}f(y)dy\quad \mbox{ for any Borel subset }A\subset\R^d
\end{equation}
and minimizing the quantity
$$\int_{\R^d}|T(x)-x|f(x)dx$$
among all the maps satisfying this condition. This means that we have a collection of particles, distributed with density $f$ on $\R^d$, that have to be moved, so that they arrange according to a new distribution, whose density is prescribed and is $g$. The movement has to be chosen so as to minimize the average displacement. The map $T$ describes the movement (that we must choose in an optimal way), and $T(x)$ represents the destination of the particle originally located at $x$. The constraint on $T$ precisely accounts for the fact that we need to reconstruct the density $g$. In  the following, we will always define, similarly to \eqref{push-forward}, the image measure of a measure $\mu$ on $X$ (measures will indeed replace the densities $f$ and $g$ in the most general formulation of the problem) through a measurable map $T:X\to Y$: it is the measure denoted by $T_\#\mu$ on $Y$ and caracterized by
\begin{gather*}
T_\#\mu(A)=\mu(T^{-1}(A)) \quad\mbox{for every measurable set } A,\\
\mbox{ or }\int_Y\phi \,\,d\left(T_\#\mu\right)=\int_X \phi\circ T\,\,d\mu\quad\mbox{for every measurable function } \phi.
\end{gather*}
The problem of Monge has stayed with no solution (does a minimizer exist? how to characterize it?\dots)  till the progress made in the 1940s. Indeed, only with the work by Kantorovich (1942) it has been inserted into a suitable framework which gave the possibility to approach it and, later, to find that solutions actually exist and to study them. The problem has been widely generalized, with very general cost functions $c(x,y)$ instead of the Euclidean distance $|x-y|$ and more general measures and spaces. For simplicity, here we will not try to present a very wide theory on generic metric spaces, manifolds and so on, but we will deal only with the Euclidean case.

\section{Primal and dual problems}

In what follows we will suppose $\Omega$ to be a (very often compact) domain of $\R^d$ and the cost function $c:\Omega\times\Omega\tto [0,+\infty[$ will be supposed continuous and symmetric (i.e. $c(x,y)=c(y,x)$).

\subsection{Kantorovich and Monge problems}

The generalization that appears as natural from the work of Kantorovich (\cite{Kantorovich}) of the problem raised by Monge is the following:

\begin{prob}
Given two probability measures $\mu$ and $\nu$ on $\Omega$ and a cost function $c:\Omega\times\Omega\tto[0,+\infty]$ we consider the problem
\begin{equation}\label{kantorovich}
(K)\quad\min\left\{\int_{\Omega\times\Omega}\!\!c\,\,d\gamma\left|\gamma\in\Pi(\mu,\nu)\right.\right\},
\end{equation}
where
$\Pi(\mu,\nu)$ is the set of the so-called {\it transport plans}, i.e. $\Pi(\mu,\nu)=\{\gamma\in\pical(\Omega\times\Omega):\,(p^+)_{\#}\gamma=\mu,\,(p^-)_{\#}\gamma=\nu,\}$
where $p^+$ and $p^-$ are the two projections of $\Omega\times\Omega$
onto $\Omega$. These probability measures over $\Omega\times\Omega$ are an alternative way to describe the displacement of the particles of $\mu$: instead of saying, for each $x$, which is the destination $T(x)$ of the particle originally located at $x$, we say for each pair $(x,y)$ how many particles go from $x$ to $y$. It is clear that this description allows for more general movements, since from a single point $x$ particles can a priori move to different destinations $y$. If multiple destinations really occur, then this movement cannot be described through a map $T$. Notice that the constraints on   $(p^\pm)_{\#}\gamma$ exactly mean that we restrict our attention to the movements that really take particles distributed according to the distribution $\mu$ and move them onto the distribution $\nu$. 

The minimizers for this problem are called {\it
optimal transport plans} between $\mu$ and $\nu$. Should $\gamma$ be of
the form $(id\times T)_{\#}\mu$ for a measurable map
$T:\Omega\tto\Omega$ (i.e. when no splitting of the mass occurs), the map $T$ would be called {\it optimal
transport map} from $\mu$ to $\nu$.
\end{prob}

\begin{rem}
It can be easily checked that if $(id\times T)_{\#}\mu$ belongs to $\Pi(\mu,\nu)$ then $T$ pushes $\mu$ onto $\nu$ (i.e. $\nu(A)=\mu(T^{-1}(A))$ for any Borel set $A$) and the functional takes the form $\int c(x,T(x))\mu(dx),$ thus generalizing Monge's problem.
\end{rem}

This generalized problem by Kantorovich is much easier to handle than the original one proposed by Monge: for instance in the Monge case we would need existence of at least a map $T$ satisfying the constraints. This is not verified when $\mu=\delta_0$, if $\nu$ is not a single Dirac mass. On the contrary, there always exist transport plan in $\Pi(\mu,\nu)$ (for instance $\mu\otimes\nu\in\Pi(\mu,\nu)$). Moreover, one can state that $(K)$ is the relaxation of the original problem by Monge: if one considers the problem in the same setting, where the competitors are transport plans, but sets the functional at $+\infty$ on all the plans that are not of the form $(id\times T)_{\#}\mu$, then one has a functional on $\Pi(\mu,\nu)$ whose relaxation is the functional in $(K)$ (see \cite{ambpra}).

Anyway, it is important to notice that an easy use of the Direct Method of Calculus of Variations (i.e. taking a minimizing sequence, saying that it is compact in some topology - here it is the weak convergence of probability measures - finding a limit, and proving semicontinuity (or continuity) of the functional we minimize, so that the limit is a minimizer) proves that a minimum does exist. 

As a consequence, if one is interested in the problem of Monge, the question may become``does this minimum come from a transport map $T$?''. Actually, if the answer to this question is yes, then it is evident that the problem of Monge has a solution, which also solves a wider problem, that of minimizing among transport plans. In some cases proving that the optimal transport plan comes from a transport map (or proving that there exists at least one optimal plan coming from a map) is equivalent to proving that the problem of Monge has a solution, since very often the infimum among transport plans and among transport maps is the same. Yet, in the presence of atoms, this is not always the case, but we will not insist any more on this degenerate case.

\subsection{Duality}

Since the problem (K) is a linear optimization under linear constraints, an important tool will be duality theory, which is typically used for convex problems. We will find a dual problem (D) for (K) and exploit the relations between dual and primal.

The first thing we will do is finding a formal dual problem, by means of an inf-sup exchange.

First express the constraint  $\gamma\in\Pi(\mu,\nu)$ in the following way : notice that, if $\gamma$ is a non-negative measure on $\Omega\times\Omega$, then we have
$$\sup_{\phi,\,\psi}\int\phi \,d\mu +\int \psi \,d\nu - \int \left(\phi(x)+\psi(y)\right) \,d\gamma = \begin{cases}0&\mbox{ if }\gamma\in\Pi(\mu,\nu)\\
+\infty&\mbox{ otherwise }\end{cases}.$$

Hence, one can remove the constraints on $\gamma$ if he adds the previous sup, since if they are satisfied nothing has been added and if they are not one gets $+\infty$ and this will be avoided by the minimization.
Hence we may look at the problem we get and interchange the inf in $\gamma$ and the sup in $\phi,\psi$:
\begin{multline*}
\min_\gamma \int c\,\,d\gamma +\sup_{\phi,\psi}\left(\!\!\int\phi \,d\mu +\int \psi \,d\nu - \int (\phi(x)+\psi(y)) \,d\gamma\!\right)=\\\sup_{\phi,\psi}\int\phi \,d\mu +\int \psi \,d\nu +\inf_\gamma \int\left(c(x,y)-(\phi(x)+\psi(y))\right) \,d\gamma.
\end{multline*}
Obviously it is not always possible to exchange inf and sup, and the main tool to do it is a theorem by Rockafellar requiring concavity in one variable, convexity in the other one, and some compactness assumption. We will not investigate anymore whether in this case these assumptions are satisfied or not. But the result is true. 

Afterwards, one can re-write the inf in $\gamma$ as a constraint on $\phi$ and $\psi$, since one has

$$\inf_{\gamma\geq 0} \int\left(c(x,y)-(\phi(x)+\psi(y))\right) \,d\gamma
= \begin{cases}0&\mbox{ if }\phi(x)+\psi(y)\leq c(x,y) \mbox{ for all }(x,y)\in\Omega\times\Omega\\
-\infty&\mbox{ otherwise }\end{cases}.$$

This leads to the following dual optimization problem.

\begin{prob}
Given the two probabilities $\mu$ and $\nu$ on $\Omega$ and the cost function $c:\Omega\times\Omega\tto[0,+\infty]$ we consider the problem
\begin{equation}\label{dual}
(D)\quad\max\left\{\int_{\Omega}\!\phi \,d\mu+\!\!\int_{\Omega}\!\psi \,d\nu\;\middle|\;\phi\in L^1(\mu),\psi\in L^1(\nu)\,:\,\phi(x)\!+\!\psi(y)\leq c(x,y) \mbox{ for all }(x,y)\in\Omega\!\times\!\Omega\right\},
\end{equation}
\end{prob}

This problem does not admit a straightforward existence result, since the class of admissible functions lacks compactness. Yet, we can better understand this problem and find existence once we have introduced the notion of $c-$transform (a kind of generalization of the well-known Legendre transform).

\begin{defi}
Given a function $\chi:\Omega\tto\overline{\R}$ we define its {\it $c-$transform} (or $c-$conjugate function) by
$$\chi^c(y)=\inf_{x\in\Omega}c(x,y)-\chi(x).$$
Moreover, we say that a function $\psi$ is {\it $c-$concave} if there exists $\chi$ such that $\psi=\chi^c$ and we denote by $\Psi_c(\Omega)$ the set of $c-$concave functions.
\end{defi}

It is quite easy to realize that, given a pair $(\phi,\psi)$ in the maximization problem (D), one can always replace it with $(\phi,\phi^c)$, and then with $(\phi^{cc},\phi^c)$, and the constraints are preserved and the integrals increased. Actually one could go on but it is possible to prove that $\phi^{ccc}=\phi^c$ for any function $\phi$. This is the same as saying that $\psi^{cc}=\psi$ for any $c-$concave function $\psi$, and this prefectly recalls what happens for the Legendre transform of convex funtions (which corresponds to the particular case $c(x,y)=x\cdot y$). 

A consequence of these considerations is the following well-known result 
\begin{prop}
We have
\begin{equation}\label{duality formula}
\min(K)=\max_{\psi\in\Psi_c(\Omega)}\int_{\Omega}\psi\,\,d\mu+\int_{\Omega}\psi^c\,\,d\nu,
\end{equation}
where the max on the right hand side is realized. In particular the minimum value of $(K)$ is a convex function of $(\mu,\nu)$, as it is a supremum of linear functionals.
\end{prop}

\begin{defi}
The functions $\psi$ realizing the maximum in \eqref{duality formula} are called {\it Kantorovich potentials} for the transport from $\mu$ to $\nu$. This is in fact a small abuse, because usually this term is used only in the case $c(x,y)=|x-y|$, but it is usually understood in the general case as well.
\end{defi}

Notice that any $c-$concave function shares the same modulus of continuity of the cost $c$. This is the reason why one can prove existence for (D) (which is the same of the right hand side problem in the previous proposition), by applying Ascoli-Arzelà's Theorem.

In, particular, in the case $c(x,y)=|x-y|^p$, if $\Omega$ is bounded with diameter $D$, any $\psi\in\Psi_c(\Omega)$ is $pD^{p-1}-$Lipschitz continuous. Notice that the case where $c$ is a power of the distance is actually of particular interest and two values of the exponent $p$ are remarkable: the cases $p=1$ and $p=2$. In these two cases we provide characterizations for the set of $c-$concave functions. Let us denote by $\Psi_{(p)}(\Omega)$ the set of $c-$concave functions with respect the cost $c(x,y)=|x-y|^p/p$. It is not difficult to check that
\begin{gather*}
\psi\in\Psi_{(1)}(\Omega) \Longleftrightarrow \psi\mbox{  is a $1$-Lipschitz function;}\\
\psi\in\Psi_{(2)}(\Omega) \Longrightarrow \;x\mapsto \frac{x^2}{2}-\psi(x)\text{  is a convex function; if $\Omega=\R^d$ this is an equivalence.}
\end{gather*}

\subsection{The case $c(x,y)=|x-y|$}

The case $c(x,y)=|x-y|$ shows a lot of interesting features, even if from the point of the existence of an optimal map $T$ it is one of the most difficult. A first interesting property is the following:
\begin{prop}
For any $1-$Lipschitz function $\psi$ we have $\psi^c=-\psi$. In particular, Formula \ref{duality formula} may be re-written as
$$\min(K)=\max(D)=\max_{\psi\in\Lip_1}\int_{\Omega}\psi\,\,d(\mu-\nu).$$
\end{prop}

The key point of the previous proposition is proving $\psi^c=-\psi$. This is easy if one considers that $\psi^c(y)=\inf_x |x-y|-\psi(x)\leq -\psi(x)$ (taking $x=y$), but also $\psi^c(y)=\inf_x |x-y|-\psi(x)\geq \inf_x |x-y|-|x-y|+\psi(y)=\psi(y)$ (making use of the Lipschitz behaviour of $\psi$).

Another peculiar feature of this case is the following:

\begin{prop}\label{minimal flow}
Consider the problem
\begin{equation}\label{bidual}
(B)\quad \min \left\{M(\lambda)\;\middle|\; \lambda\in\M^d(\Omega);\; \nabla\cdot\lambda=\mu-\nu\right\},
\end{equation}
where $M(\lambda)$ denotes the mass of the vector measure $\lambda$ and the divergence condition is to be read in the weak sense, with Neumann boundary conditions, i.e. $-\int\nabla\phi\!\cdot\! d\lambda=\int\phi \,d(\mu-\nu)$ for any $\phi\in C^1(\overline\Omega)$. If $\Omega$ is convex then it holds
$$\min (K) = \min (B).$$
\end{prop}

This proposition links the Monge-Kantorovich problem to a minimal flow problem which has been first proposed by Beckmann in \cite{beck}, under the name of {\it continuous transportation model}. He did not know this link, as Kantorovich's theory was being developed independently almost in the same years. In Section 2.1 we will see some details more on this model and on the possibility of generalizing it to the case of distances $c(x,y)$ coming from Riemannian metrics. In particular, in the case of a nonconvex $\Omega$, $(B)$ would be equivalent to a Monge-Kantorovich problem where $c$ is the geodesic distance on $\Omega$.

To have an idea of why these equivalences between (B) and (K) hold true, one can look at the following considerations.

First, a formal computation.
We take the problem (B) and re-write the constraint on $\lambda$ by means of the quantity
$$\sup_{\phi}\int-\nabla\phi\cdot \,\,d\lambda  +\int \phi \,d(\mu-\nu)  = \begin{cases}0&\mbox{ if }\nabla\cdot\lambda=\mu-\nu\\
+\infty&\mbox{ otherwise }\end{cases}.$$

Hence one can write (B) as
\begin{equation*}
\min_\lambda M(\lambda)+\sup_\phi\int-\nabla\phi\cdot \,\,d\lambda  +\int \phi \,d(\mu-\nu)=\sup_{\phi} \int\phi \,d(\mu-\nu)+\inf_\lambda M(\lambda)-\int\nabla\phi\cdot \,\,d\lambda,
\end{equation*}
where inf and sup have been exchanged formally as in the previous computations. After that one notices that 
$$\inf_\lambda M(\lambda)-\int\nabla\phi\cdot \,\,d\lambda=\inf_\lambda\int d|\lambda|\left(1-\nabla\phi\cdot\frac{\,d\lambda}{d|\lambda|}\right)= \begin{cases}0&\mbox{ if }|\nabla\phi|\leq 1\\
-\infty&\mbox{ otherwise }\end{cases}$$
and this leads to the dual formulation for (B) which gives
$$\sup_{\phi\,:\,|\nabla\phi|\leq 1}\int_{\Omega}\phi\,\,d(\mu-\nu).$$
Since this problem is exactly the same as (D) (a consequence of the fact that $\Lip_1$ functions are exactly those functions whose gradient is smaller than $1$), this gives the equivalence between (B) and (K).

Most of the considerations above, especially those on the problem (B) do not hold for costs other than the distance $|x-y|$. The only possible generalizations I know concern either a cost $c$ which comes from a Riemannian distance $k(x)$ (i.e. $c(x,y)=\inf \{\int_0^1 k(\sigma(t))|\sigma'(t)|dt\,:\,\sigma(0)=x,\sigma(1)=y\}$, which gives a problem (B) with $\int k(x)d|\lambda|$ instead of $M(\lambda)$) or the fact that $p-$homogeneous costs may become $1-$homogeneous through the introduction of time as an extra variable (see \cite{PhDChloe}). Some more details on the problem (B) can be found in the lectures notes on ``Models and applications of optimal transport in economics, traffic and urban planning'' of this same Summer School, \cite{GrenobleModels}.

\subsection{$c(x,y)=h(x-y)$ with $h$ strictly convex and the existence of an optimal $T$} 

We summarize here some useful results for the case where the cost $c$ is of the form $c(x,y)=h(x-y)$, for a strictly convex function $h$.

The main tool is the duality result. If we have equality between the minimum of (K) and the maximum of (D) and both extremal values are realized, one can consider an optimal transport plan $\gamma$ and a Kantorovich potential $\psi$ and write
$$\psi(x)+\psi^c(y)\leq c(x,y) \mbox{ on }\Omega\times\Omega \mbox{ and }\psi(x)+\psi^c(y)= c(x,y) \mbox{ on }\spt\gamma.$$
The equality on $\spt\gamma$ is a consequence of the inequality which is valid everywhere and of 
$$\int c\,\,d\gamma=\int\psi \,\,d\mu+\int\psi^c \,\,d\nu=\int(\psi(x)+\psi^c(y))\,d\gamma,$$
which implies equality $\gamma-$a.e. These functions being continuous, the equality passes to the support of the measure.

Once we have that, let us fix a point $(x_0,y_0)\in\spt\gamma$. One may deduce from the previous computations that
$$x\mapsto \psi(x)-h(x-y_0)\quad\mbox{ is minimal at }x=x_0$$
and, if $\psi$ is differentiable at $x_0$, one gets $\nabla\psi(x_0)\in\partial h(x_0-y_0).$
For a strictly convex function $h$ one may inverse the relation passing to $\nabla h^*$ thus getting
$$x_0-y_0=\nabla h^*(x_0)=(\partial h)^{-1}(x_0).$$
This solves several questions concerning the transport problem with this cost, provided $\psi$ is differentiable a.e. with respect to $\mu$. This is usually guaranteed by requiring $\mu$ to be absolutely continuous with respect to the Lebsgue measure, and using the fact that $\psi$ may be proven to be Lipschitz. Then, one may use the previous computation to deduce that, for every $x_0$, the point $y_0$ such that $(x_0,y_0)\in\spt\gamma$ is unique (i.e. $\gamma$ is of the form $(id\times T)_\#\mu$ where $T(x_0)=y_0$). Moreover, this also gives uniqueness of the optimal trasport plan and of the gradient of the Kantorovich potential.

We may summarize everything in the following theorem:
\begin{teo}\label{maintransport}
Given $\mu$ and $\nu$ probability measures on a domain $\Omega\subset\R^d$ there exists an optimal transport plan $\pi$. It is unique and of the form $(id\times T)_{\#}\mu$, provided $\mu$ is absolutely continuous. Moreover there exists also at least a Kantorovich potential $\psi$, and the gradient $\nabla\psi$ is uniquely determined $\mu-$a.e. (in particular $\psi$ is unique up to additive constants, provided the density of $\mu$ is positive a.e. on $\Omega$). The optimal transport map $T$ and the potential $\psi$ are linked by $T(x)=x-(\nabla h^*)(\nabla\psi(x))$. Moreover we have $\psi(x)+\psi^c(T(x))=c(x,T(x))$ for $\mu-$a.e. $x$. Conversely, every map $T$ which is of the form $T(x)=x-(\nabla h^*)(\nabla\psi(x))$ for a function $\psi\in\Psi_c(\Omega)$ is an optimal transport plan from $\mu$ to $T_{\#}\mu$.
\end{teo}
\begin{rem}
Actually, the existence of an optimal transport map is true under weaker assumptions: we can replace the condition of being absolutely continuous with the condition ``$\mu(A)=0$ for any $A\subset\R^d$ such that $\haus^{d-1}(A)<+\infty$'' or with any condition which ensures that the non-differentiability set of $\psi$ is negligible. In the theorem we used the Lipschitz behavior of $\psi\in\Psi_c$ and applied Rademacher Theorem, but $c-$concave functions are often more regular than only Lipschitz.
\end{rem}
\begin{rem}
In Theorem \ref{maintransport} only the part concerning the optimal map $T$ is not symmetric in $\mu$ and $\nu$: hence the uniqueness of the Kantorovich potential is true even if it $\nu$ (and not $\mu$) has positive density a.e. (since one can retrieve $\psi$ from $\psi^c$ and viceversa).
\end{rem}
\begin{rem}\label{casep2}
Theorem \ref{maintransport} may be particularized to the quadratic case $c(x,y)=|x-y|^2/2$, thus getting the existence of an optimal transport map 
$$T(x)=x-\nabla\psi(x)=\nabla\left(\frac{x^2}{2}-\psi(x)\right)=\nabla\phi(x)$$ 
for a convex function $\phi$. By using the converse implication (sufficient optimality conditions), this also proves the existence and the uniqueness of a gradient of a convex function transporting $\mu$ onto $\nu$. This well known fact has been investigated first by Brenier (see \cite{Brenier polar}) and is often known as Brenier's Theorem. 

Let us moreover notice that a specific approach for the case $|x-y|^2$, based on the fact that we can withdraw the parts of the cost depending on $x$ or $y$ only and maximize $\int \!x\!\cdot\! y \,d\gamma$, gives the same result in a easier way: we actually get $\phi(x_0)+\phi^*(y_0)=x_0\cdot y_0$ for a convex function $\phi$ and its Legendre transform $\phi^*$ and we deduce $y_0\in\partial\phi(x_0)$. 
\end{rem}

All the costs of the form $c(x,y)=|x-y|^p$ with $p>1$ fall under Theorem \ref{maintransport}. 

We finish the part dedicated to positive results  by noticing that the same method may not be used if $h$ is only convex, or at least does not give results as strong as what it does if $h$ is strictly convex. Yet, there is anyway something which is known 
for the case $c(x,y)=|x-y|$.  The results are a bit weaker (and much harder) and are summarized below (this is the classical Monge case and we refer to \cite{ambpra}, even if several different proofs have been provided by different methods). Notice that a lot of literature is currently being dedicated to the case of other norms than the Euclidean one and other distance functions. 
\begin{teo}\label{maintransportp1}
Given $\mu$ and $\nu$ probability measures on a domain $\Omega\subset\R^d$ there exists at least an optimal transport plan $\pi$ for the cost $c(x,y)=|x-y|$. Moreover, one of such plans is of the form $(id\times T)_{\#}\mu$ provided $\mu$ is absolutely continuous. There exists a Kantorovich potential $\psi$, and its gradient is unique $\mu-$a.e.and we have $\psi(x)-\psi(T(x))=|x-T(x)|$ for $\mu-$a.e. $x$, for any choice of optimal $T$ and $\psi$.
\end{teo}
Here the absolute continuity assumption is essential to have existence of an optimal transport map, in the sense that in general it cannot be replaced by weaker assumptions as in the strictly convex case.

Morevoer, we can provide a counter-exemple showing that in general it is necessary that $\mu$ does not give mass to ``small'' sets.

\begin{ex}
Set
$$\mu=\huno\res A\;\mbox{ and }\;\nu=\frac{\huno\res B+\huno\res C}{2}$$
where $A$, $B$ and $C$ are three vertical parallel segments in $\R^2$ whose vertexes lie on the two line $y=0$ and $y=1$ and the abscissas are $0$, $1$ and $-1$, respectively, and $\huno$ is the $1-$dimensional Haudorff measure. It is clear that no transport plan may realize a cost better than $1$ since, horizontally, every point needs to be displaced of a distance $1$. Moreover, one can get a sequence of maps $T_n: A\to B\cup C$ by dividing $A$ into $2n$ equal segments $(A_i)_{i=1,\dots,2n}$ and $B$ and $C$ into $n$ segments each, $(B_i)_{i=1,\dots,n}$ and $(C_i)_{i=1,\dots,n}$ (all ordered downwards). Then define $T_n$ as a piecewise affine map which sends $A_{2i-1}$ onto $B_i$ and $A_{2i}$ onto $C_i$. In this way the cost of the map $T_n$ is less than $1+1/n$, which implies that the infimum of the Kantorovich problem is $1$, as well as the infimum on transport maps only. Yet, no map $T$ may obtain a cost $1$, as this would imply that all points are sent horizontally,but this cannot respect the push-forward constraint. On the other hand, the transport plan associated to $T_n$ weakly converge to the transport plan $\frac 12 T^+_{\#}\mu+\frac 12 T^-_{\#}\mu,$ where $T^{\pm}(x)=x\pm e$ and $e=(1,0)$. This transport plan turns out to be the only optimal transport plan and its cost is $1$.
\end{ex}
Notice that the same construction provides also an example of the relaxation procedure leading from Monge to Kantorovich.

  \begin{center}
  \setlength{\unitlength}{0.25cm}
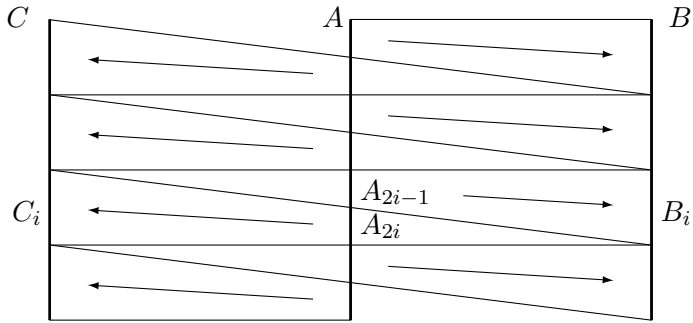
\begin{picture}(35,18)(-15,0)
\put(2,14.875){\vector(16,-1){12}}
\put(0,16){\line(16,0){16}}
\put(0,14){\line(16,-2){16}}
\put(0,14){\line(-16,2){16}}
\put(-2,13.125){\vector(-16,1){12}}
\put(0,12){\line(-16,0){16}}
\put(2,10.875){\vector(16,-1){12}}
\put(0,12){\line(16,0){16}}
\put(0,10){\line(16,-2){16}}

\put(0,10){\line(-16,2){16}}
\put(-2,9.125){\vector(-16,1){12}}
\put(0,8){\line(-16,0){16}}

\put(0,6){\line(-16,2){16}}
\put(-2,5.125){\vector(-16,1){12}}
\put(0,4){\line(-16,0){16}}

\put(0,2){\line(-16,2){16}}
\put(-2,1.125){\vector(-16,1){12}}
\put(0,0){\line(-16,0){16}}

\put(6,6.625){\vector(16,-1){8}}
\put(0,8){\line(16,0){16}}
\put(0,6){\line(16,-2){16}}

\put(2,2.875){\vector(16,-1){12}}
\put(0,4){\line(16,0){16}}
\put(0,2){\line(16,-2){16}}
\put(-1.5,15.6){$A$}
\put(16.9,15.6){$B$}
\put(-18.3,15.6){$C$}

\put(16.5,5.3){$B_i$}
\put(-18,5.3){$C_i$}

\put(0.5,6.4){$A_{2i-1}$}
\put(0.5,4.6){$A_{2i}$}
\thicklines
\put(0,0){\line(0,1){16}}
\put(16,0){\line(0,1){16}}
\put(-16,0){\line(0,1){16}}
\put(0,0){\line(0,1){16}}
\put(16,0){\line(0,1){16}}
\put(-16,0){\line(0,1){16}}
\put(0,0){\line(0,1){16}}
\put(16,0){\line(0,1){16}}
\put(-16,0){\line(0,1){16}}
\end{picture}
\end{center}

\section{Wasserstein distances and spaces}

Starting from the values of the problem $(K)$ in \eqref{kantorovich} we can define a set of distances over $\pical(\Omega)$. For any $p\geq 1$ we can define
$$W_p(\mu,\nu)=\big(\min(K)\text{ with }c(x,y)=|x-y|^p\big)^{1/p}.$$
We recall that, by Duality Formula, we have
\begin{equation}\label{duality formula2}
\frac{1}{p}W_p^p(\mu,\nu)= \sup_{\psi\in\Psi_{(p)}(\Omega)} \int_{\Omega}\psi\,\,d\nu + \int_{\Omega}\psi^c\,\,d\mu.
\end{equation}
\begin{teo}
If $\Omega$ is compact, for any $p\geq 1$ the function $W_p$ is in fact a distance over $\pical(\Omega)$ and the convergence with respect to this distance is equivalent to the weak convergence of probability measures. In particular any functional $\mu\mapsto W_p(\mu,\nu)$ is continuous with respect to weak topology.
\end{teo}

To prove that the convergence according to $W_p$ is equivalent to weak convergence one first establish this result for $p=1$, through the use of the duality with the functions in $\Lip_1$. Then it is possible to use the inequalities between the distances $W_p$ (see below) to extend the result to a general $p$.

The case of a noncompact $\Omega$ is a little more difficult. First, the distance must be defined only on a subset of the whole space of probability measures, to avoid infinite values. We will use the space of probabilities with finite $p-$th momentum:
$$\mathcal{W}_p(\Omega)=\{\mu\in\pical(\Omega)\,:\,m_p(\mu):=\int_{\Omega}|x|^p\mu(dx)<+\infty\}.$$
\begin{teo}\label{conWp}
For any $p\geq 1$ the function $W_p$ is a distance over $\mathcal{W}_p(\Omega)$ and, given a measure $\mu$ and a sequence $(\mu_n)_n$ in $\W_p(\Omega)$, the following are equivalent:
\begin{itemize}
\item $\mu_n\to\mu$ according to $W_p$;
\item $\mu_n\deb\mu$ and $m_p(\mu_n)\to m_p(\mu)$;
\item $\int_{\Omega}\phi\,\,d\mu_n\to\int_{\Omega}\phi\,\,d\mu$ for any $\phi\in C^0(\Omega)$ whose growth is at most of order $p$ (i.e. there exist constants $A$ and $B$ depending on $\phi$ such that $|\phi(x)|\leq A+B|x|^p$ for any $x$).
\end{itemize}
\end{teo}

Notice that, as a consequence of H\"older (or Jensen) inequalities, the Wasserstein distances are always ordered, i.e. $W_{p_1}\leq W_{p_2}$ if $p_1\leq p_2$. Reversed inequalities are possible only if $\Omega$ is bounded, and in this case we have, if set $D=\diam(\Omega)$, for $p_1\leq p_2$,
$$W_{p_1}\leq W_{p_2}\leq D^{1-p_1/p_2}W_{p_1}^{p_1/p_2}.$$

From the monotone behavior of Wasserstein distances with respect to $p$ it is natural to introduce the distance $W_{\infty}$: set $\W_{\infty}(\Omega)=\{\mu\in\pical(\Omega)\,:\,\spt(\mu)\,\mbox{ is bounded }\}$ (obviously if $\Omega$ itself is bounded one has $\W_{\infty}(\Omega)=\pical(\Omega)$) and then
$$W_{\infty}(\mu,\nu)=\inf\left\{\gamma-\mathrm{esssup}_{x,y\in\Omega\times\Omega}|x-y|\,:\,\gamma\in\Pi(\mu,\nu)\right\}.$$
Here $\gamma-\mathrm{esssup}$ denotes the essential sup with respect to $\gamma$, i.e. the norm in the space $L^\infty(\Omega\!\times\Omega;\gamma)$, which is the same, for continuous functions such as $|x-y|$, as the maximal value on the support of $\gamma$.
It is easy to check that $W_p\nearrow W_{\infty}$ and it is interesting to study the metric space $\W_{\infty}(\Omega)$. Curiously enough, this supremal problem in optimal transport theory, even if quite natural, has not deserved much attention, up to the very recent paper \cite{ChaDePJut}.

The $W_{\infty}$ convergence is stronger than any $W_p$ convergence and hence also than the weak convergence of probability measures. The converse is not true and the convergence in $W_{\infty}$ turns out to be actually rare: consequently there is a great lack of compactness in $\W_{\infty}$. For instance it is not difficult to check that, if we set $\mu_t=t\delta_{x_0}+(1-t)\delta_{x_1}$, where $x_0\neq x_1\in\Omega$, we have $W_{\infty}(\mu_t,\mu_s)=|x_0-x_1|$ if $t\neq s$. This implies that the balls $B(\mu_t,|x_0-x_1|/2)$ are infinitely many disjoint balls in $\W_{\infty}$ and prevents compactness.

The following statement summarizes the compactness properties of the spaces $W_p$ for $1\leq p\leq \infty$ and its proof is a direct application of the considerations above and of Theorem \ref{conWp}.
\begin{prop}\label{compactWp}
For $1\leq p <\infty$ the space $W_p(\Omega)$ is compact if and only if $\Omega$ itself is compact. Moreover, for an unbounded $\Omega$ the space $\W_p(\Omega)$ is not even locally compact. The space $W_{\infty}(\Omega)$ is neither compact nor locally compact for any choice of $\Omega$ with $\#\Omega>1$.
\end{prop}

\section{Geodesics, continuity equation and displacement convexity}

\subsection{Metric derivatives in Wasserstein spaces}

We are concerned in this sections with several properties linked to the curves in the Wasserstein space $W_p$. For this subject the main reference is \cite{AmbGigSav}. Before giving the main result we are interested in, we recall the definition of metric derivative, which is a concept that may be useful when studying curves which are valued in generic metric spaces.

\begin{defi}
Given a metric space $(X,d)$ and a curve $\gamma:[0,1]\to X$ we define {\it metric derivative} of the curve $\gamma$ at time $t$ the quantity
\begin{equation}\label{metric derivative}
|\gamma'|(t) = \lim_{s \to t} \frac{d(\gamma(s),\gamma(t))}{|s - t|},
\end{equation}
provided the limit exists.
\end{defi}
As a consequence of Rademacher Theorem it can be seen (see \cite{AmbTil}) that for any Lipschitz curve the metric derivative exists at almost every point $t\in[0,1]$. We will be concerned quite often with metric derivatives of curves which are valued in the space $\W_p(\Omega)$.

\begin{defi}
If we are given a Lipschitz curve $\mu:[0,1]\tto \W_p(\Omega)$, we define velocity field of the curve
any vector field $v:[0,1]\times\Omega\tto \R^d$ such that for a.e. $t\in[0,1]$ the vector field
$v_t=v(t,\cdot)$ belongs to $[L^p(\mu_t)]^d$ and the {\it continuity equation}
$$\frac{d}{dt}\mu_t+\nabla\cdot (v\cdot\mu_t)=0$$
is satisfied in the sense of distributions: this means that for all
$\phi\in C^1_c(\Omega)$ and any $t_1<t_2\in[0,1]$ it holds
$$\io \phi\, \,d\mu_{t_2}-\io\phi \,\,d\mu_{t_1}=\int_{t_1}^{t_2}ds\io \nabla\phi\cdot v_s\, \,d\mu_s,$$
or, equivalently, in differential form:
$$\frac{d}{d t}\io \phi \,\,d\mu_{t}=\io \nabla\phi\cdot v_t\, \,d\mu_t\qquad\mbox{ for a.e. }t\in[0,1].$$
We say that $v$ is the \emph{tangent} field to the curve $\mu_t$ if, for a.e. $t$,
$v_t$ has minimal $[L^p(\mu_t)]^d$ norm for any $t$ among all the velocity fields (actually this is not the true definition of a tangent vector field, since this would involve the definition of a tangent space for the ``manifold'' $\W_p$, but it is in this case the same).
\end{defi}

The following proposition is concerned with the existence of tangent fields and comes from Theorem~8.3.1 and Proposition~8.4.5 in \cite{AmbGigSav}.

\begin{teo}\label{tangent field}
If $p>1$ and $\mu=(\mu_t)_t$ is a curve in $\Lip([0,1];W_p(\Omega))$ then there exists a unique vector field $v$ characterized by
\begin{gather}\label{continuity equation}
\frac{\partial }{\partial t}\mu+\nabla\cdot(v\cdot\mu)=0,\\
||v_t||_{L^p(\mu_t)}\leq |\mu'|(t)\,\mbox{ for a.e. }t,
\end{gather}
where the continuity equation is satisfied in the sense of distributions as previously explained. Moreover, if \eqref{continuity equation} holds for a family of vector fields $(v_t)_t$ with $||v_t||_{L^p(\mu_t)}\leq C$ then $\mu\in\Lip([0,1];\W_p(\Omega))$ and $|\mu'|(t)\leq ||v_t||_{L^p(\mu_t)}$ for a.e. $t$.
\end{teo}

To have an idea of the meaning of the previous theorem and of the relationship between curves of measures and the continuity equation some considerations could be useful.

Actually, at least when the vector fields $v_t$ are regular enough, the solution of the continuity equation $\partial\mu/\partial t +\nabla\cdot(v\cdot\mu)=0,$ are obtained by taking the images of the initial measure $\mu_0$ through the maps $\sigma(t,\cdot)$ obtained by taking the solution of
$$\begin{cases}\sigma'(t,x)=v_t(\sigma(t,x)),\\
			\sigma(0,x)=x.\end{cases}$$
This explains why the vector field $v_t$ is called ``velocity field'' of the curve $\mu_t$: if every particle follows at each time $t$ the velocity field $v_t$, then the position of all the particles at time $t$ reconstructs exactly the measure $\mu_t$ that appears in the continuity equation together with $v_t$ !

Think for a while to the case of two time steps only: there are two measures $\mu_t$ and $\mu_{t+h}$ and there are several ways for moving the particles so as to reconstruct the latter from the former. It is exactly as when we look for a transport. One of these trasnports is optimal in the sense that it minimizes $\int |T(x)-x|^p\mu_t(dx)$ and the value of this integral equals $W_p^p(\mu_t,\mu_{t+h})$. If we call $v_t(x)$ the ``discrete velocity of the particle located at $x$ at time $t$, i.e. $v_t(x)=(T(x)-x)/h$, one has $||v_t||_{L^p(\mu_t)}=\frac 1 h W_p(\mu_t,\mu_{t+h})$.  The result of the previous theorem may be easily guessed as obtainable as a limit as $h\to 0$.

\subsection{Geodesics and geodesic convexity}

Once we know about curves in their generality, it is interesting to think about geodesics. The following result is a characterization of geodesics in $W_p(\Omega)$ when $\Omega$ is a convex domain in $\R^d$. This procedure is also known as {\it McCann's linear interpolation}.

\begin{teo}\label{geodesics in Wp}
All the spaces $\W_p(\Omega)$ are length spaces and if $\mu$ and $\nu$ belong to $\W_p(\Omega)$, and $\gamma$ is an optimal transport plan from $\mu$ to $\nu$ for the cost $c_p(x,y)=|x-y|^p$, then the curve
$$\mu^{\gamma}(s)=(p_s)_{\#}\gamma$$
where $p_s:\Omega\times\Omega\to\Omega$ is given by $p_s(x,y)=x+s(y-x)$, is a constant-speed geodesic  from $\mu$ to $\nu$. In the case $p>1$ all the constant-speed geodesics are of this form, and if $\mu$ is absolutely continuous, then there is only one geodesic and it has the form
$$
\mu(s)= [(1 - s)id + s T]_{\#} \mu,
$$
where $T$ is the optimal transport map from $\mu$ to $\nu$.
\end{teo}

By means of this characterization of geodesics we can also define the useful concept of displacement convexity introduced by McCann in \cite{MC}.
\begin{defi}\label{displacement_interpolation}
Given a functional $F:\W_p(\Omega)\cap L^1\tto [0,+\infty]$, we say that it is {\it displacement convex} if all the maps $t \mapsto F(\mu^{\gamma}(t))$ are convex on $[0,1]$ for every choice of $\mu$ and $\nu$ in $\W_p(\Omega)$ and $\gamma$ optimal transport plan from $\mu$ to $\nu$ with respect to $c(x,y)=|x-y|^p$.
\end{defi}

The following well-known result provides a wide set of displacement convex functionals. In the case $p=2$ this result is due to McCann (\cite{MC}), while the generalization to any $p$ can be found in \cite{AmbGigSav}.

\begin{teo}\label{displ conv result}
Consider the following functionals on the space $\W_p(\Omega)$, where $\Omega$ is any convex subset of $\R^N$:
\begin{eqnarray*}
J^1(\mu)&=&\begin{cases}\int_{\Omega} f(u(x))\, d x & \text{if } \mu=u\cdot\lcal^d\\
                     + \infty & \text{ if $\mu$ is not absolutely continuous};\end{cases}\\
J^2(\mu)&=& \int_{\Omega} V(x)\,\mu(dx);\\
J^3(\mu)&=& \int_{\Omega}\int_{\Omega} w(x-y)\mu(dx)\mu(dy).
\end{eqnarray*}

Suppose that $f : [0,+\infty] \to [0,+\infty]$ is a convex and superlinear lower semicontinuous function with $f(0) = 0$, and that $V:\Omega\tto[0,+\infty]$ and $w:\R^d\tto[0,+\infty]$ are convex functions. Then the functionals $J^2$ and $J^3$ are displacement convex in $\W_p(\Omega)$ and the functional $J^1$ is displacement convex provided the map
$$
r \mapsto r^d f(r^{-d})
$$
is convex and non-increasing on $]0,+ \infty[$.
\end{teo}

\section{Monge-Ampère equation and regularity}

The final issue that we'll approach in these lecture notes will be concerned with some regularity properties of $T$ and $\psi$ (the optimal transport map and the Kantorovich potential, respectively) and their relations with the densities of $\mu$ and $\nu$. We will consider only the quadratic case $c(x,y)=|x-y|^2/2$, because it is the one where more results have been proven. Very recent results for generic costs have been developed by Ma, Trudinger, Wang, Loeper, Figalli\dots They require some very rigid assumptions on the costs, so that, surprinsingly enough, the quadratic cost is one of the few power that satisfies the suitable hypotheses.

It is easy -- just by a change-of-variables formula -- to transform the equality $\nu=T_{\#}\mu$ into the PDE $v(T(x))=u(x)/|JT|(x)$, where $u$ and $v$ are the densities of $\mu$ and $\nu$ (which have to be supposed regular enough) and $J$ denotes the determinant of the Jacobian matrix. Recalling that we may write $T=\nabla\phi$ with $\phi$ convex (Remark \ref{casep2}), we get the Monge-Ampère equation
\begin{equation}\label{mongeampere}
M\phi=\frac{u}{v(\nabla\phi)},
\end{equation}
where $M$ denotes the determinant of the Hessian
$$M\phi=\det H\phi= \det \left[\frac{\partial^2\phi}{\partial x_i\,\partial x_j}\right]_{i,j} .$$
This equation up to now is satisfied by $\phi=\frac{x^2}{2}-\psi$ in a formal way only. We define various notions of solutions for \eqref{mongeampere}:
\begin{itemize}
\item we say that $\phi$ satisfies \eqref{mongeampere} in the Brenier sense if $(\nabla\phi)_{\#}u\cdot\lcal^d=v\cdot\lcal^d$ (and this is actually the sense to be given to this equation);
\item we say that $\phi$ satisfies \eqref{mongeampere} in the Alexandroff sense if $H\phi$, which is always a positive measure for $\phi$ convex, is absolutely continuous and its density satisfies \eqref{mongeampere} a.e.;
\item we say that $\phi$ satisfies \eqref{mongeampere} in the viscosity sense if it satisfies the usual comparison properties required by viscosity theory but restricting the comparisons to regular convex test functions (since $M$ is in fact monotone just when restricted to positively definite matrices);
\item we say that $\phi$ satisfies \eqref{mongeampere} in the classical sense if it is of class $C^2$ and the equation holds pointwise.
\end{itemize}
Notice that any notion except the first may be also applied to the more general equation $M\phi=f$, while the first one just applies to this specific transportation case.
The results we want to use are well summarized in Theorem 50 of \cite{villani}:
\begin{teo}
If $u$ and $v$ are $C^{0,\alpha}(\Omega)$ and are both bounded from above and from below on the whole $\Omega$ by positive constants and $\Omega$ is a convex open set, then the unique Brenier solution $\phi$ of \eqref{mongeampere} belongs to $C^{2,\alpha}(\Omega)\cap C^{1,\alpha}(\overline{\Omega})$ and $\phi$ satisfies the equation in the classical sense (hence also in the Alexandroff and viscosity senses).
\end{teo}
Even if this precise statement is taken from \cite{villani}, we just detail a possible bibliographical path to arrive at this result. It is not easy to deal with Brenier solutions, so the idea is to consider viscosity solutions, for which it is in general easy to prove existence by Perron's method. Then prove some regularity result on viscosity solutions, up to getting a classical solution. After that, once we have a classical convex solution to Monge-Ampère equation, this will be a Brenier solution too. Since this is unique (up to additive constants) we have got a regularity statement for Brenier solutions.
We can find results on viscosity solutions in \cite{caf1}, \cite{caf2} and \cite{caf3}. In \cite{caf1} some conditions to ensure strict convexity of the solution of $M\phi=f$ when $f$ is bounded from above and below are given. In \cite{caf2} for the same equation it is proved $C^{1,\alpha}$ regularity provided we have strict convexity. In this way the term $u/v(\nabla\phi)$ becomes a $C^{0,\alpha}$ function and in \cite{caf3} it is proved $C^{2,\alpha}$ regularity for solutions of $M\phi=f$ with $f\in C^{0,\alpha}$.

\end{document}